\documentclass[12pt]{amsart}
\usepackage{amsmath, amssymb, amsthm, bm, mathtools, enumitem, caption,microtype,mathrsfs}
\usepackage{hyperref}

\usepackage[noabbrev,capitalize]{cleveref}
\usepackage{adjustbox}
\crefname{equation}{}{}
\usepackage{ytableau}
\usepackage{graphicx, tikz}
\usepackage[textheight=8.6in, textwidth=6.5in]{geometry}

\usepackage{microtype}

\newtheorem{theorem}{Theorem}[section]
\newtheorem{lemma}[theorem]{Lemma}

\newtheorem{proposition}[theorem]{Proposition}

\newtheorem*{conjecture*}{Conjecture}
\newtheorem*{problem}{Problem}

\theoremstyle{definition}
\newtheorem{definition}[theorem]{Definition}

\theoremstyle{remark}
\newtheorem*{remark}{Remark}
\newtheorem*{remarks}{Remarks}
\newtheorem*{example}{Example}

\numberwithin{equation}{section}

\DeclarePairedDelimiter\abs{\lvert}{\rvert}

\DeclareMathOperator{\Tr}{Tr}

\newcommand{\HH}{\mathbb H}

\newcommand{\C}{\mathbb C}
\newcommand{\SL}{\mathrm{SL}}
\newcommand{\F}{\mathbb F}

\newcommand{\Q}{\mathbb Q}
\newcommand{\RR}{\mathbb R}
\newcommand{\Z}{\mathbb Z}

\renewcommand{\abs}[1]{\left\vert #1 \right \vert}

\newcommand{\spt}{\operatorname{spt}}

\newcommand{\Supersing}{\operatorname{SS}}

\newcommand{\e}{\operatorname{e}}
\newcommand{\Mcal}{\mathcal M}
\newcommand{\Qcal}{\mathcal Q}

\newcommand{\Dmod}{\Delta_{\mathrm{mod}}}
\newcommand{\ch}{\chi_{12}}
\newcommand{\leg}[2]{\left(\frac{#1}{#2}\right)}
\newcommand{\norm}[1]{\lVert #1\rVert}
\newcommand{\vech}{\boldsymbol e}

\title[Smallest parts, CM traces, and congruences]
{The smallest-parts function and elliptic curves}
\author{Ajit Singh}
\address{Department of Mathematics $\&$ Computing, Indian Institute of Technology (Indian School of Mines) Dhanbad, India}
\email{ajit94@iitism.ac.in}

\begin{document}
\begin{abstract} 
Recently, Ono established a new connection between partition congruences and the arithmetic of elliptic curves by expressing the partition function \(p(n)\) as a CM trace on the modular curve \(X_0(6)\). This approach yields new proofs of Ramanujan-type congruences for the partition function modulo powers of \(5\), \(7\), and \(11\), revealing how the geometry of supersingular elliptic curves and the arithmetic of oriented optimal embeddings govern these congruence properties. Motivated by Ono's work, we develop an analogous passage from CM traces to congruences
modulo \(5\), \(7\), and \(13\) for the for Andrews's smallest-parts partition function \(\spt(n)\), employing the established theta-lift and modular-trace methods suited to this setting. With \(P\) the Bruinier--Ono weight-zero
weak Maass function and \(f\) an explicit weakly holomorphic modular
function on \(\Gamma_0(6)\), the Ahlgren--Andersen formula expresses
\(\spt(n)\) as one twelfth of the discriminant \(1-24n\) CM trace of
\(f-P\).  We explain the theta-lift normalization, the contribution of
the two orientations, and the harmonic-Maass-form uniqueness argument.
Following Garvan's work, we construct integral forms \(G_\ell\) of weight
\(\ell+1\) on \(\Gamma_0(\ell)\), compute their Fricke transformations,
and prove the exact identities
\[
\mathcal T_5(G_5)=-5E_6,\qquad
\mathcal T_7(G_7)=-7E_8,\qquad
\mathcal T_{13}(G_{13})=13E_{14}.
\]
Their reductions yield the three congruences on the Andrews's
smallest-parts partition function.  We also
give a supersingular specialization of the trace of \(f\), with actual
reduction multiplicities, and explain the role of imprimitive forms,
the conductor restriction in CM reduction, and the distinction between
the primes \(11\) and \(13\). These results have been verified in Lean using AxiomProver.
 \end{abstract}
\maketitle
\section{Introduction and Statement of Results}
For a positive integer \(n\), let \(p(n)\) denote the number of partitions
of \(n\), and let \(\spt(n)\) denote the total number of occurrences of
the smallest parts in all partitions of \(n\).  For example, the partitions of \(4\) are
$
4, 3+1, 2+2, 2+1+1, 1+1+1+1,
$
and their smallest-part multiplicities are \(1,1,2,2,4\), respectively.  Hence
\(\spt(4)=10\). The generating function for \(\spt(n)\) is given by
\begin{align}
	S(q):=\sum_{m\geq 0}\spt(m)q^m&=\sum_{m\ge1}
	\frac{q^m}{(1-q^m)^2(q^{m+1};q)_\infty}
	=\frac1{(q;q)_\infty}
	\sum_{m\ge1}\frac{q^m(q;q)_{m-1}}{1-q^m}.\label{eq:correctGF2}
\end{align}
where we have used  the {\it $q$-Pochhammer symbol}
\begin{equation}
	(a;q)_n:= \begin{cases} 1   \ \ \ \ \ & {\text {\rm if}}\ n=0,  \\
		(1-a)(1-aq)\cdots (1-aq^{n-1}) \ \ \ \ \ &{\text {\rm if}}\ n\in \Z_+,\\
		\prod_{j=0}^{\infty}(1-aq^j) \ \ \ \ \ &{\text {\rm if}}\ n=+\infty.
	\end{cases}
\end{equation}
Throughout, we put \(p(0):=1\), \(\spt(0):=0\), and extend both functions by zero at negative integer arguments. Extensive research has been devoted to the study of the $\spt(n)$-function; see, for example, \cite{AA, ADR,AGL, AGL1,FolsomOno,Garvan,Ono, Rolen}.

Recently, Ono \cite{OnoEC} prove that for each $n\geq 1$, the partition function $p(n)$ as a CM trace on $X_0(6)$ of
the discriminant $\delta_n = 1 -24n$ invariants of a weight 0 weak Maass function $P$ that records
where CM elliptic curves sit on $X(1)$, together with their canonical first-order “CM tangent”, the diagonal local slope of the CM isogeny relation on $X(1)\times X(1)$. This approach yields new proofs of Ramanujan-type congruences for the partition function modulo powers of \(5\), \(7\), and \(11\), revealing how the geometry of supersingular elliptic curves and the arithmetic of oriented optimal embeddings govern these congruence properties.  In \cite{Andrews}, Andrews introduced \(\spt(n)\) and proved three congruences modulo
\(5\), \(7\), and \(13\).  The generating function $S(q)$ is related
to a harmonic Maass form of weight \(3/2\), through Bringmann's work on
partition statistics \cite{Bringmann}.  This places smallest parts in
the same broad arithmetic setting as the ordinary partition function,
although the relevant modular weights and the congruence primes differ.

The finite algebraic formula of Bruinier and Ono \cite{BO} writes
\(p(n)\) as a CM trace of a weight-zero weak Maass function.
Ahlgren and Andersen \cite{AA} obtained the corresponding formula for
\(\spt(n)\). Motivated by Ono's recent work \cite{OnoEC}, we develop an analogous passage from CM traces to congruences for the smallest-parts function, employing the established theta-lift and modular-trace methods suited to this setting.
\subsection{Notation and the two level-six functions}
For positive integers $k$ and \(\tau\in\HH\), the weight $2k$ Eisenstein series (see Ch. 1 of \cite{CBMS}) is
\begin{equation}\label{Eisenstein}
	E_{2k}(\tau) :=1-\frac{4k}{B_{2k}}\sum_{n=1}^{\infty}\sigma_{2k-1}(n)q^n
\end{equation}
where $B_{2k}$ is the  $2k$-th Bernoulli number, $\sigma_{\nu}(n):=\sum_{d\mid n}d^{\nu},$ and $q:=e^{2\pi i \tau}.$ The first few examples are
\[
E_2(\tau)=1-24\sum_{m\ge1}\sigma_1(m)q^m,\quad
E_4(\tau)=1+240\sum_{m\ge1}\sigma_3(m)q^m,\quad
E_6(\tau)=1-504\sum_{m\ge1}\sigma_5(m)q^m.
\]
Apart from $E_2,$ each $E_{2k}$ is a weight $2k$ holomorphic modular form on $\SL_2(\Z).$  
We also require the Dedekind eta function
\begin{align*}
	\eta(\tau) := q^{\frac{1}{24}}\prod_{n\geq 1}(1-q^n).
\end{align*}
Put
\begin{equation}\label{eq:H6}
	H_6(\tau):=\bigl(\eta(\tau)\eta(2\tau)\eta(3\tau)\eta(6\tau)\bigr)^2
\end{equation}
and define
\begin{equation}\label{eq:Fminus2}
	F(\tau):=
	\frac{E_2(\tau)-2E_2(2\tau)-3E_2(3\tau)+6E_2(6\tau)}
	{2H_6(\tau)}.
\end{equation}
The Bruinier--Ono function is
\begin{equation}\label{eq:P}
	P(\tau)=-
	\left(\frac{1}{2\pi i}\frac{d}{d\tau}+\frac{1}{2\pi y}\right)F(\tau),\qquad y=Im(\tau).
\end{equation}
The additional modular function for smallest parts is
\begin{equation}\label{eq:f}
	f(\tau)=
	\frac{E_4(\tau)-4E_4(2\tau)-9E_4(3\tau)+36E_4(6\tau)}
	{24H_6(\tau)}.
\end{equation}
The first few terms are
\[
F(\tau)=q^{-1}-10-29q-104q^2+\cdots,\qquad
f(\tau)=q^{-1}+12+77q+376q^2+\cdots.
\]

For \(n\ge1\), let \(D_n:=1-24n\), and let
\begin{equation}\label{eq:Qset}
	\Qcal_{D_n}^{(1)}
	:=
	\left\{[a,b,c]\in\Z^3:
	a>0,\ b^2-4ac=D_n,\ 6\mid a,\ b\equiv1\pmod{12}\right\}.
\end{equation}
All forms in this set are positive definite.  We include
\emph{all integral forms}, without a primitivity restriction.
For \(Q=[a,b,c]\), write
\begin{equation}\label{eq:CMpoint}
	\alpha_Q:=\frac{-b+\sqrt{D_n}}{2a}\in\HH,\qquad
	\Gamma=\Gamma_0(6)/\{\pm I\}.
\end{equation}
For a \(\Gamma_0(6)\)-invariant function \(u\), define
\begin{equation}\label{eq:CMtrace}
	\Tr_{D_n}(u):=\sum_{Q\in\Gamma\backslash\Qcal_{D_n}^{(1)}}u(\alpha_Q).
\end{equation}
These are finite sums.  The notation \(\Tr_{D_n}\) always means this CM
trace; a trace of modular forms over group cosets will be denoted by
\(\mathcal T_\ell\).

\begin{theorem}\cite[Ahlgren--Andersen]{AA}\label{thm:CMmain}
	For every integer \(n\ge1\),
	\begin{equation}\label{eq:sptCM}
		\spt(n)=\frac1{12}\Tr_{D_n}(f-P)
		=\frac1{12}
		\sum_{Q\in\Gamma\backslash\Qcal_{D_n}^{(1)}}
		\bigl(f(\alpha_Q)-P(\alpha_Q)\bigr).
	\end{equation}
	Equivalently,
	\begin{equation}\label{eq:ftrace}
		\Tr_{D_n}(f)=12\spt(n)+(24n-1)p(n).
	\end{equation}
	The CM values in \eqref{eq:sptCM} are algebraic.  They belong to the
	ring class fields of the relevant CM orders, and the values in each
	content packet are permuted by the corresponding CM class-group action.
\end{theorem}

Theorem~\ref{thm:CMmain} is the
smallest-parts counterpart of the Bruinier--Ono formula
\begin{equation}\label{eq:partitionCM}
	\Tr_{D_n}(P)=(24n-1)p(n).
\end{equation}
The proof in Section~\ref{sec:CMproof} explains the factors \(2\) and
\(12\), which are important when comparing different normalizations.

\begin{theorem}\cite[Andrews]{Andrews}\label{thm:congruences}
	For every integer \(n\ge0\),
	\begin{align}
		\spt(5n+4)&\equiv0\pmod5,\label{eq:spt5}\\
		\spt(7n+5)&\equiv0\pmod7,\label{eq:spt7}\\
		\spt(13n+6)&\equiv0\pmod{13}.\label{eq:spt13}
	\end{align}
\end{theorem}

We prove Theorem~\ref{thm:congruences} uniformly by Garvan's
construction \cite{Garvan}.  It produces
\(G_\ell\in M_{\ell+1}(\Gamma_0(\ell))\cap\Z[[q]]\), for primes
\(\ell\ge5\), whose full-group traces satisfy
\begin{equation}\label{eq:introK}
	\mathcal T_5(G_5)=-5E_6,\qquad
	\mathcal T_7(G_7)=-7E_8,\qquad
	\mathcal T_{13}(G_{13})=13E_{14}.
\end{equation}
Here \(E_8=E_4^2\) and \(E_{14}=E_4^2E_6\) and $M_{k}(\Gamma_0(N))$ is the space of modular forms of weight $k$ and level $N$.  The Fricke transformation
shows that \(\mathcal T_\ell(G_\ell)\equiv G_\ell\pmod\ell\), and the
Fourier expansion of \(G_\ell\) extracts the required spt progression.

\subsection{Scope and organization}
The main CM formula and the three congruences are established
theorems.  Our purpose is an exposition in which their definitions,
normalizations, and proof inputs are stated together.  The non-elementary
inputs are Bringmann's harmonic completion, the Bruinier--Funke
theta-lift theorem, the Bruinier--Ono partition trace, and standard
facts about complex multiplication and integral modular curves.
The specialized coefficient calculations, the uniqueness argument,
the Fricke calculation, and the deduction of the congruences are
included below.\\
\noindent
The organization of the paper is as follows. In Section~\ref{sec:generating}, we recall preliminary notions and results on the harmonic completion.  Section~\ref{sec:CMproof} proves \eqref{eq:ftrace}.  Sections~\ref{sec:atkin} and \ref{sec:congruenceproof} construct the modular trace and prove Theorem~\ref{thm:congruences}.  Section~\ref{sec:geometry} proves the integral regularity of \(f\) and
gives the supersingular specialization with explicit hypotheses, and
Section~\ref{sec:comparison} compares the analytic and arithmetic
features of \(\mathrm{spt}(n)\) with the partition setting.
\section*{Acknowledgements}
\noindent The author thanks Ken Ono for suggesting this problem. The author also gratefully acknowledges financial support from the INSPIRE Faculty Research Grant (IFA-24MA204) and the Anusandhan National Research Foundation (ANRF) Early Career Research Grant (ECRG/2025/009350/PMS). The author used Google Gemini to improve the exposition of the paper. All suggestions were carefully reviewed by the author, who takes full responsibility for the content of this paper. After completing the first draft, the author shared the manuscript with Ken Ono, who facilitated its formal verification using AxiomProver. The Lean formalization and verification of the proofs were carried out entirely by the AxiomProver system.
\section{Nuts and Bolts}\label{sec:generating}
In this section, we recall some definitions and useful results on harmonic completion (for example, see \cite{Bringmann, BFOR}). We also  fix some notations as well. We define for $n\ge0$
\begin{equation}\label{eq:an}
	a(n):=12\spt(n)+(24n-1)p(n).
\end{equation}
Thus \(a(0)=-1\).  At nonintegral or negative arguments we use the
convention $a(t)=0$ if $t\notin\Z_{\ge0}$. By taking the logarithmic derivative of \(\eta\), we obtain that
\[
\frac1{2\pi i}\frac{\eta'(\tau)}{\eta(\tau)}=\frac1{24}E_2(\tau),
\]
and therefore
\begin{equation}\label{eq:E2eta}
	-\frac{E_2(\tau)}{\eta(\tau)}
	=\sum_{n\ge0}(24n-1)p(n)q^{n-1/24}.
\end{equation}
It follows that
\begin{equation}\label{eq:alpha}
	\alpha(\tau)
	:=12q^{-1/24}S(q)-\frac{E_2(\tau)}{\eta(\tau)}
	=\sum_{n\ge0}a(n)q^{n-1/24}.
\end{equation}
In particular, principal part of $\alpha$ is $-q^{-1/24}$.

We use the real Dirichlet character
\begin{equation}\label{eq:chi12}
	\ch(m)=\leg{12}{m}
	:=
	\begin{cases}
		1,&m\equiv1,11\pmod{12},\\
		-1,&m\equiv5,7\pmod{12},\\
		0,&\gcd(m,6)>1.
	\end{cases}
\end{equation}
Euler's pentagonal theorem gives the useful unary theta expansion
\begin{equation}\label{eq:etaUnary}
	\eta(\tau)=\sum_{m\ge1}\ch(m)q^{m^2/24}.
\end{equation}

 Throughout, we let $\tau = x + i y \in \HH$, where $x,y\in \RR$. We recall that for $k\in\mathbb{R}$, the weight $k$ hyperbolic Laplacian operator on
$\mathbb{H}$ is defined by
\[
\Delta_k:=-y^2\left(\frac{\partial^2}{\partial x^2}
+\frac{\partial^2}{\partial y^2}\right)
+iky\left(\frac{\partial}{\partial x}
+i\frac{\partial}{\partial y}\right).
\]
Harmonic Maass forms of weight \(k\) are real-analytic functions on the upper half-plane \(\mathbb{H}\) that lie in the kernel of the weight-\(k\) hyperbolic Laplacian \(\Delta_k\). The transformation law for half-integral weights requires additional factors, which we now introduce. Let \(\left(\frac{c}{d}\right)\) denote the extended Legendre symbol, and let the square root be taken on the principal branch of the holomorphic square-root function. For each odd integer \(d\), define
\[
\varepsilon_d :=
\begin{cases}
	1, & d\equiv 1\pmod{4},\\
	i, & d\equiv 3\pmod{4}.
\end{cases}
\tag{4.1}
\]
With the above notations, we have the following definition.
\begin{definition}
If $k\in\frac12\mathbb{Z}$, then a weight $k$ harmonic Maass form on
a subgroup $\Gamma=\Gamma_0(N)$ for some $N\in\mathbb{N}$, where
$4\mid N$ if $k\in\frac12\mathbb{Z}\setminus\mathbb{Z}$ throughout,
is any smooth function $f:\mathbb{H}\to\mathbb{C}$ satisfying the
following properties:

\begin{enumerate}
	\item[(i)] For all
	$\begin{pmatrix}a&b\\c&d\end{pmatrix}\in\Gamma$ and all
	$\tau\in\mathbb{H}$, we have
	\[
	f\left(\frac{a\tau+b}{c\tau+d}\right)
	=
	\begin{cases}
		(c\tau+d)^k f(\tau), & k\in\mathbb{Z},\\[4pt]
		\left(\dfrac{c}{d}\right)
		\varepsilon_d^{-2k}(c\tau+d)^k f(\tau),
		& k\in\frac12+\mathbb{Z}.
	\end{cases}
	\]
	
	\item[(ii)] We have
	\[
	\Delta_k(f)=0.
	\]
	
	\item[(iii)] There exists a polynomial
	$P_f(\tau)\in\mathbb{C}[q^{-1}]$ such that
	\[
	f(\tau)-P_f(\tau)=O\left(e^{-\varepsilon v}\right)
	\]
	as $v\to\infty$ for some $\varepsilon>0$. Analogous conditions are
	required at all cusps.
\end{enumerate}
\end{definition}
\noindent
Note that in particular weakly holomorphic modular forms are harmonic Maass forms. A weakly holomorphic modular form is holomorphic on $\HH$ and meromorphic at the cusps.
\begin{theorem}\cite[Theorem 1.1]{Bringmann}\label{thm:completion}
	Let
	\begin{equation}\label{eq:Mminus}
		\Mcal^-(\tau)=\frac3{\sqrt\pi}\sum_{m\ge1}
		\ch(m)m\,
		\Gamma\left(-\frac12,\frac{\pi m^2y}{6}\right)q^{-m^2/24},
	\end{equation}
	where \(\Gamma(s,x)=\int_x^\infty t^{s-1}e^{-t}\,dt\).
	Then
	\begin{equation}\label{eq:M}
		\Mcal(\tau)=\alpha(\tau)+\Mcal^-(\tau)
	\end{equation}
	is a harmonic Maass form of weight \(3/2\) with inverse eta multiplier:
	\begin{equation}\label{eq:Mtransform}
		\Mcal(\gamma\tau)=
		\nu_\eta(\gamma)^{-1}(c\tau+d)^{3/2}\Mcal(\tau),
		\qquad
		\gamma=\begin{pmatrix}a&b\\c&d\end{pmatrix}\in\SL_2(\Z).
	\end{equation}
	Here \(\eta(\gamma\tau)=\nu_\eta(\gamma)(c\tau+d)^{1/2}\eta(\tau)\),
	with the usual compatible choices of square roots.
	Moreover,
	\begin{equation}\label{eq:shadow}
		\xi_{3/2}\Mcal=-\frac{3\sqrt6}{\pi}\eta,\quad\text{where}~~
		\xi_kG=2iy^k\,\overline{\frac{\partial G}{\partial\bar\tau}}.
	\end{equation}
\end{theorem}

\begin{remarks} 
(1) We note that the modularity assertion is the harmonic-completion result obtained
	from Bringmann's work \cite{Bringmann}; the normalization used here is
	recorded in \cite[Theorem~2.1]{Garvan} and
	\cite[(3.1)--(3.4)]{AA}.\\
(2)	The constants can be checked from the equivalent completion
	\begin{equation}\label{eq:EichlerIntegral}
		\Mcal^-(\tau)=\frac{6\sqrt{3i}}{\pi}
		\int_{-\bar\tau}^{i\infty}\frac{\eta(w)}{(w+\tau)^{3/2}}\,dw.
	\end{equation}
	Insert \eqref{eq:etaUnary}, write \(w=-x+it\), and then put \(u=t+y\).
	For \(b=m^2/24\), the integral becomes
	\[
	\frac{6\sqrt3}{\pi}q^{-b}
	\int_{2y}^{\infty}u^{-3/2}e^{-2\pi bu}\,du
	=\frac{3m}{\sqrt\pi}
	\Gamma\left(-\frac12,\frac{\pi m^2y}{6}\right)q^{-m^2/24}.
	\]
	This establishes \eqref{eq:Mminus}.  \\
(3) By differentiating the incomplete gamma function termwise yields \eqref{eq:shadow}.
	In particular, the two generator transformations are
	\begin{equation}\label{eq:MST}
		\Mcal(\tau+1)=e^{-2\pi i/24}\Mcal(\tau),\qquad
		\Mcal(-1/\tau)=e^{\pi i/4}\tau^{3/2}\Mcal(\tau).
	\end{equation}
\end{remarks}

\begin{proposition}\label{lem:unique}
	Let \(H\) be a harmonic Maass form of weight \(3/2\) satisfying
	\eqref{eq:Mtransform}, with holomorphic shadow at the cusp.  If its
	holomorphic part has zero principal part, then \(H=0\).
\end{proposition}
\begin{proof}
	Let \(h := \xi_{3/2}H\), where \(\xi_k(F) := 2iy^k\overline{\frac{\partial F}{\partial\bar\tau}}\) is the antilinear differential operator of Bruinier and Funke \cite[Proposition~3.2]{BFtwo}. By the mapping properties of \(\xi_{3/2}\), since \(H\) has weight \(3/2\) and multiplier \(\nu_\eta^{-1}\), the shadow \(h\) transforms as a modular form of dual weight \(2 - 3/2 = 1/2\) on \(\SL_2(\Z)\) with the conjugate multiplier \(\overline{\nu_\eta^{-1}} = \nu_\eta\).
	Moreover, because \(H\) is harmonic, it is annihilated by the weight \(3/2\) hyperbolic Laplacian \(\Delta_{3/2} = -\xi_{1/2}\xi_{3/2}\) \cite[Equation~(3.1)]{BFtwo}, so \(\xi_{1/2}h = 0\), which means \(h\) is holomorphic on \(\HH\). Under the translation matrix \(T = \left(\begin{smallmatrix} 1 & 1 \\ 0 & 1 \end{smallmatrix}\right)\), the eta multiplier takes the value \(\nu_\eta(T) = e^{\pi i/12} = e^{2\pi i/24}\) \cite[Chapter~3]{Knoff}. Thus \(h(\tau+1) = e^{2\pi i/24}h(\tau)\), and the Fourier expansion of \(h\) is supported on exponents in \(1/24 + \Z\):
	\[
	h(\tau) = \sum_{n \in \Z} b(n) q^{n + 1/24}.
	\]
	By hypothesis, \(h\) is holomorphic at the cusp \(\infty\), so \(b(n) = 0\) for all \(n < 0\). Consequently, \(h(\tau) = \sum_{n \ge 0} b(n)q^{n+1/24} = O(q^{1/24})\) as \(y = \operatorname{Im}(\tau) \to \infty\), so \(h\) vanishes at \(\infty\). Since \(\eta(\tau)\) also transforms with weight \(1/2\) and multiplier \(\nu_\eta\) and is nonvanishing on \(\HH\), the quotient \(h/\eta\) is a holomorphic modular function of weight \(0\) on \(\SL_2(\Z)\) satisfying \(h(\tau)/\eta(\tau) = b(0) + O(q)\) as \(y \to \infty\). Because any holomorphic modular function of weight \(0\) on \(\SL_2(\Z)\) that is bounded at \(\infty\) is constant \cite[Chapter~VII, Proposition~1]{Serre}, there exists a constant \(c \in \C\) such that \(h = c\eta\). In particular, \(h\) is a cusp form of weight \(1/2\).\\
	\noindent We next show that \(h = 0\) by computing its Petersson norm \(\norm{h}_{\mathrm{Pet}}^2\) via the Bruinier--Funke pairing \cite[Proposition~3.5]{BFtwo} and Stokes' theorem. For \(T > 1\), let $\mathcal{F}_T$ 
	be the standard truncated fundamental domain for \(\SL_2(\Z)\). Using the definition \(\overline{\xi_{3/2}H} = -2iy^{3/2}\frac{\partial H}{\partial\bar\tau}\) together with the holomorphy of \(h\) and the relation \(d\bar\tau \wedge d\tau = 2i\,dx \wedge dy\), we have the exact differential form identity
	\begin{align*}
	|h(\tau)|^2 y^{1/2} \frac{dx\,dy}{y^2}
	&= h(\tau)\,\overline{\xi_{3/2}H(\tau)}\,y^{-3/2}\,dx\,dy
	= -2i\,h(\tau)\frac{\partial H}{\partial\bar\tau}\,dx\,dy\\
	&= -\frac{\partial(Hh)}{\partial\bar\tau}\,d\bar\tau \wedge d\tau
	= -d\bigl(H(\tau)h(\tau)\,d\tau\bigr).
\end{align*}
	Since \(H\) transforms with weight \(3/2\) and multiplier \(\nu_\eta^{-1}\), \(h\) transforms with weight \(1/2\) and multiplier \(\nu_\eta\), and \(d(\gamma\tau) = (c\tau+d)^{-2}\,d\tau\), the \(1\)-form \(\omega := H(\tau)h(\tau)\,d\tau\) is strictly \(\SL_2(\Z)\)-invariant, i.e., 
	\[
	H(\gamma\tau)h(\gamma\tau)\,d(\gamma\tau)
	= \nu_\eta(\gamma)^{-1}(c\tau+d)^{3/2} H(\tau) \cdot \nu_\eta(\gamma)(c\tau+d)^{1/2} h(\tau) \cdot (c\tau+d)^{-2}\,d\tau
	= H(\tau)h(\tau)\,d\tau.
	\]
	Moreover, \(\omega\) is smooth throughout \(\HH\). Applying Stokes' theorem on \(\mathcal{F}_T\) yields
	\[
	\int_{\mathcal{F}_T} |h(\tau)|^2 y^{1/2} \frac{dx\,dy}{y^2} = -\int_{\partial\mathcal{F}_T} H(\tau)h(\tau)\,d\tau.
	\]
	By the \(\SL_2(\Z)\)-invariance of \(H(\tau)h(\tau)\,d\tau\), the line integrals along these paired boundary components cancel completely, leaving only the integral along the top horizontal segment \(y = T\) (oriented from \(x = 1/2\) to \(x = -1/2\)):
	\[
	\int_{\mathcal{F}_T} |h(\tau)|^2 y^{1/2} \frac{dx\,dy}{y^2} = \int_{-1/2}^{1/2} H(x+iT)h(x+iT)\,dx.
	\]
	To evaluate the limit as \(T \to \infty\), we decompose \(H = H^+ + H^-\) into its holomorphic and nonholomorphic parts \cite[Equations~(3.2)--(3.4)]{BFtwo}. Because \(H(\tau+1)= e^{-2\pi i/24}H(\tau)\), the holomorphic part \(H^+\) has a Fourier expansion supported on exponents in \(-1/24 + \Z\):
	\[
	H^+(\tau) = \sum_{n \gg -\infty} c^+(n) q^{n - 1/24}.
	\]
	Since \(H^+\) has zero principal part at \(\infty\), we have \(c^+(n) = 0\) for all \(n \le 0\), and hence \(H^+(\tau) = \sum_{n \ge 1} c^+(n)q^{n-1/24} = O(q^{23/24}) = O(e^{-2\pi(23/24)T})\) as \(T \to \infty\). Meanwhile, since \(h = c\eta = O(q^{1/24})\) is a cusp form, the nonholomorphic part \(H^-\) is given by an incomplete Gamma series of the form \(\sum_{m \ge 1} c^-(m)\Gamma(-1/2, 4\pi (m/24)y)q^{-m/24}\) \cite[Equation~(3.3)]{BFtwo}. Using the asymptotic \(\Gamma(s, Y) \sim Y^{s-1}e^{-Y}\) as \(Y \to \infty\), each term in \(H^-(x+iT)\) is bounded by \(O\bigl(T^{-3/2}e^{-4\pi(m/24)T}\cdot e^{2\pi(m/24)T}\bigr) = O\bigl(T^{-3/2}e^{-2\pi(m/24)T}\bigr)\), so \(H^-(x+iT)\) decays exponentially as \(T \to \infty\). Together with \(h(x+iT) = O(e^{-2\pi T/24})\), it follows that
	\[
	\norm{h}_{\mathrm{Pet}}^2 = \lim_{T \to \infty} \int_{\mathcal{F}_T} |h(\tau)|^2 y^{1/2} \frac{dx\,dy}{y^2} = \lim_{T \to \infty} \int_{-1/2}^{1/2} H(x+iT)h(x+iT)\,dx = 0.
	\]
	Because the integrand \(|h(\tau)|^2 y^{-3/2}\) is continuous and nonnegative on \(\HH\), we conclude that \(h = \xi_{3/2}H = 0\). Consequently, \(\frac{\partial H}{\partial\bar\tau} = 0\), which means \(H^- = 0\) and \(H = H^+\) is a holomorphic modular form on \(\HH\).
	
	Finally, consider the product \(F(\tau) := H(\tau)\eta(\tau)\). Since \(H\) is holomorphic of weight \(3/2\) with multiplier \(\nu_\eta^{-1}\) and \(\eta\) is holomorphic of weight \(1/2\) with multiplier \(\nu_\eta\), their product \(F\) is a holomorphic modular form of weight 2 on \(\SL_2(\Z)\) with trivial multiplier (\(\nu_\eta^{-1}\nu_\eta = 1\)). At the cusp \(\infty\), we have
	$
	F(\tau) = H^+(\tau)\eta(\tau) = O(q^{23/24}) \cdot O(q^{1/24}) = O(q),
	$
	so \(F\) is holomorphic at \(\infty\). Thus \(F \in M_2(\SL_2(\Z))\). However, it is a standard consequence of the valence formula for \(\SL_2(\Z)\) \cite[Chapter~VII, Theorem~3 and Corollary~1]{Serre} that \(M_2(\SL_2(\Z)) = \{0\}\): for any nonzero \(f \in M_k(\SL_2(\Z))\).
	 Since \(\eta(\tau) \neq 0\) for all \(\tau \in \HH\), it follows that \(H = 0\).
\end{proof}

\section{CM theta lifts and the finite algebraic formula}
\label{sec:CMproof}

\subsection{Modularity and Atkin--Lehner symmetries}
Using Theorems 1.64 and 1.65 in \cite{CBMS}, we find that \(H_6\) is a weight-four modular
form with trivial character on \(\Gamma_0(6)\). Also, factors of \(H_6\) have no zeros on \(\HH\).
The numerator of \(f\) is a weight-four modular form on the same
group, so \(f\) is weakly holomorphic of weight zero.
Its Fourier expansion is integral: its numerator divided by \(24\)
has constant term \(1\) and integral remaining coefficients, while
\(H_6=q\) times a unit in \(\Z[[q]]\). We use the determinant-one representatives
\begin{align*}
	W_2=\frac1{\sqrt2}\begin{pmatrix}2&-1\\6&-2\end{pmatrix},\qquad
	W_3=\frac1{\sqrt3}\begin{pmatrix}3&1\\6&3\end{pmatrix},\qquad
	W_6=\frac1{\sqrt6}\begin{pmatrix}0&-1\\6&0\end{pmatrix}.
\end{align*}
Applying the transformations of \(\eta\) and \(E_4\) to
\eqref{eq:f} gives
\begin{equation}\label{eq:ALf}
	f|_0W_2=-f,\qquad f|_0W_3=-f,\qquad f|_0W_6=f.
\end{equation}
The same signs hold for \(P\); they follow first for \(F\), and then
from the compatibility of the Maass raising operator with the slash
action.  Explicitly, \eqref{eq:ALf} exchanges the \(d\)-isogeny terms
in the numerator, with \(d=1,2,3,6\); the signs are the character
of the Atkin--Lehner group taking values \(-1,-1,+1\) at
\(W_2,W_3,W_6\).
\subsection{The lattice and its Weil representation}
\label{subsec:lattice-weil}
We begin by introducing the rational quadratic space underlying the
theta-lift construction. For more details on the general theory of lattices, discriminant forms, and their
Weil representations, see Borcherds~\cite{Borcherds1998} and
Bruinier--Funke~\cite{BFtwo}. Let
\begin{equation}\label{eq:lattice}
	L=
	\left\{
	\begin{pmatrix}
		b & c/6\\
		-a & -b
	\end{pmatrix}
	:a,b,c\in\mathbb Z
	\right\},
	\qquad
	Q_L(X)=6\det X.
\end{equation}
For
\[
X=
\begin{pmatrix}
	b & c/6\\
	-a & -b
\end{pmatrix},\quad \text{we have}\quad Q_L(X)=ac-6b^2.
\]
The associated bilinear form is
\[
(X,Y)=Q_L(X+Y)-Q_L(X)-Q_L(Y)
=-6\operatorname{tr}(XY).
\]
Consequently, \(L\) is an even integral lattice of rank three, and
the associated rational quadratic space has signature \((1,2)\).

The dual lattice is defined by
\[
L'=\{X\in L\otimes_{\mathbb Z}\mathbb Q:
(X,Y)\in\mathbb Z
\text{ for every }Y\in L\}.
\]
A direct calculation gives
\begin{equation}\label{eq:dual-lattice}
	L'=
	\left\{
	\begin{pmatrix}
		b/12 & c/6\\
		-a & -b/12
	\end{pmatrix}
	:a,b,c\in\mathbb Z
	\right\}.
\end{equation}
It follows that the discriminant group
\[
D_L=L'/L
\]
is cyclic of order \(12\). Under the identification
\(D_L\simeq\mathbb Z/12\mathbb Z\), its finite quadratic form is
\begin{equation}\label{eq:discriminant-form}
	q_L(h)\equiv-\frac{h^2}{24}\pmod{\mathbb Z}.
\end{equation}
 The group \(\Gamma_0(6)\) acts on \(L\) by conjugation,
\(X\mapsto\gamma X\gamma^{-1}\). This action preserves the quadratic
form and induces the trivial action on \(D_L\). The finite quadratic
module \(D_L\) therefore determines a Weil representation
\[
\rho_L:\operatorname{Mp}_2(\mathbb Z)
\longrightarrow\operatorname{GL}(\mathbb C[D_L]),
\]
where \(\operatorname{Mp}_2(\mathbb Z)\) denotes the metaplectic
double cover of \(\operatorname{SL}_2(\mathbb Z)\).

Let \(\{\mathfrak e_h:h\in\mathbb Z/12\mathbb Z\}\) be the standard
basis of \(\mathbb C[D_L]\), and put \(\e(z)=e^{2\pi iz}\).
With the chosen convention for the Weil representation, its
standard generators act by
\begin{equation}\label{eq:Weil}
	\rho_L(T)\mathfrak e_h
	=\e\left(-\frac{h^2}{24}\right)\mathfrak e_h,
\end{equation}
and
\begin{equation}\label{eq:Weil-S}
	\rho_L(S)\mathfrak e_h
	=\frac{e^{\pi i/4}}{\sqrt{12}}
	\sum_{h'\bmod12}
	\e\left(\frac{hh'}{12}\right)\mathfrak e_{h'}.
\end{equation}
These formulas are the finite Fourier-transform description of the
Weil representation associated with the discriminant form
\eqref{eq:discriminant-form}; see~\cite{Borcherds1998,BFtwo}.

Let
\[
\chi_{12}(h)=\left(\frac{12}{h}\right)
\]
be the Kronecker character, and define
\begin{equation}\label{eq:character-vector}
	\mathfrak v_\chi:
	=\sum_{h\bmod12}\chi_{12}(h)\mathfrak e_h.
\end{equation}
The quadratic Gauss-sum identity
\begin{equation}\label{eq:chiGauss12}
	\sum_{h\bmod12}\chi_{12}(h)
	\e\left(\frac{hh'}{12}\right)
	=\sqrt{12}\,\chi_{12}(h')
\end{equation}
shows that \(\mathbb C\mathfrak v_\chi\) is invariant under
\(\rho_L(S)\). Since \(h^2\equiv1\pmod{24}\) whenever
\(\chi_{12}(h)\ne0\), it is also invariant under \(\rho_L(T)\).
Thus \(\mathfrak v_\chi\) spans a one-dimensional subrepresentation
of the Weil representation.
In particular, a vector-valued modular form of the form
\[
\mathbf F(\tau)=G(\tau)\mathfrak v_\chi
\]
transforms according to \(\rho_L\) if the scalar-valued function
\(G\) has the compatible inverse multiplier. This is the mechanism
we use to relate scalar-valued modular forms to the vector-valued
modular forms arising in theta-lift constructions; for more details see
\cite{BFtwo}.

\subsection{The theta-lift input and its specialization}
The Bruinier--Funke theta-lift theorem \cite[Theorem~4.5 and
Proposition~4.7]{BFtraces} associates to a weakly holomorphic
weight-zero modular function \(u\) the harmonic Maass form
\[
I(\tau,u)=\int_{\Gamma\backslash\HH}u(z)\Theta_L(\tau,z)
=\sum_{h\bmod12}I_h(\tau,u)\vech_h,
\]
with the standard theta-lift interpretation of the integral.
Here \(\Theta_L\) is the Kudla--Millson theta kernel for \(L\).
For \(m>0\), we define
\[
L_{h,m}:=\{X\in L+h:Q_L(X)=m\}.
\]
The coefficient formula is
\begin{equation}\label{eq:BFpositive}
	[q^m]I_h(\tau,u)^+
	=
	\sum_{X\in\Gamma\backslash L_{h,m}}
	\frac{u(z_X)}{\abs{\Gamma_X}},
\end{equation}
where \(z_X\) is the point corresponding to the positive real line
\(\mathbb R X\).  The negative coefficients are cusp contributions
described below.  These general theta-lift facts are the second
foundational input of the proof our main result.

The lift is equivariant under Atkin--Lehner operators:
\begin{equation}\label{eq:Iequiv}
	I_{Wh}(\tau,u)=I_h(\tau,u|_0W).
\end{equation}
On \(L'/L\), \(W_2,W_3,W_6\) act by multiplication by \(7,5,11\),
respectively.  Hence, equations \eqref{eq:ALf} and \eqref{eq:Iequiv} give
\begin{equation}\label{eq:Ichi}
	I_h(\tau,f)=\ch(h)I_1(\tau,f).
\end{equation}
For nonunit residues \(h\), either \(W_2\) or \(W_3\) fixes \(h\);
the corresponding minus sign forces \(I_h=0\).  For unit residues,
the identity follows by moving from \(h=1\).  By
\eqref{eq:Weil}--\eqref{eq:chiGauss12}, \(I_1\) has inverse eta
multiplier.
Our next lemma gives the relationship between the Fourier coefficients of the positive part of the theta lift and traces of CM values.
\begin{lemma}\label{lem:factor2}
	For \(n\ge1\),
	\[
	[q^{n-1/24}]I_1(\tau,f)^+=2\Tr_{D_n}(f).
	\]
\end{lemma}
\begin{proof}
	Every vector in \(L_{1,n-1/24}\) has the form
	\[
	X=\begin{pmatrix}b+1/12&c/6\\-a&-b-1/12\end{pmatrix}.
	\]
	Its norm equation is
	\[
	ac-\frac{(12b+1)^2}{24}=n-\frac1{24}.
	\]
	On the subset \(a>0\), the equivariant correspondence is
	\begin{equation}\label{eq:XQ}
		X\longleftrightarrow Q_X=[6a,12b+1,c].
	\end{equation}
	The equation is exactly \(\operatorname{disc}(Q_X)=D_n\), and
	\(z_X=\alpha_{Q_X}\).  There is no primitivity restriction in this
	correspondence.
	
	It is important to handle the second orientation correctly. Observe that the negation alone sends component \(1\) to component \(-1\).
	Instead, the involution
	\[
	X\longmapsto -W_6XW_6^{-1}
	\]
	preserves component \(1\) and exchanges \(a>0\) with \(a<0\).
	If \(B=b+1/12\), its value is
	\[
	\begin{pmatrix}B&-a/6\\c&-B\end{pmatrix},
	\]
	whose new parameter \(a\) is \(-c\); positive norm implies \(ac>0\).
	The two corresponding values of \(f\) agree by \(f|_0W_6=f\).
	Finally, since \(\Gamma_0(6)/\{\pm I\}\) is torsion-free, so
	\(\abs{\Gamma_X}=1\).  Formula \eqref{eq:BFpositive} now gives twice
	the trace over \eqref{eq:XQ}. This completes the proof of the lemma.
\end{proof}

\begin{lemma}\label{lem:principal}
	The principal part of \(I_1(\tau,f)\) is \(-2q^{-1/24}\).
\end{lemma}
\begin{proof}
	Here we show the cusp calculation to fix the normalization.
	The cusps for the congruence subgroup $\Gamma_0(6)$ are \(\infty,0,1/2,1/3\), with widths \(1,6,3,2\), respectively.
	We choose scaling matrices
	\[
	\sigma_\infty=I,\quad
	\sigma_0=\begin{pmatrix}0&-1\\1&0\end{pmatrix},\quad
	\sigma_{1/2}=\begin{pmatrix}1&1\\2&3\end{pmatrix},\quad
	\sigma_{1/3}=\begin{pmatrix}1&-1\\3&-2\end{pmatrix}.
	\]
	We define \(f(\sigma_{\mathfrak a}\tau):
	=\sum_r c_{\mathfrak a}(r)q^r\).
	The following table shows the cusp data for this lattice and these scalings.
	\begin{center}
		\begin{tabular}{|p{2.5cm}|p{1.0cm}|p{1.0cm}|p{1.2cm}|p{1.2cm}|}
		\hline
		$\mathfrak a$&$\infty$&0&1/2&1/3\\
		 \hline
		$w_{\mathfrak a}$&1&6&3&2\\
		\hline
		$\beta_{\mathfrak a}$&1/6&1&1/2&1/3\\
		\hline
		$w_{\mathfrak a}/\beta_{\mathfrak a}$&6&6&6&6\\
		\hline
		\text{polar term of } $f|_0\sigma_{\mathfrak a}$
		&$q^{-1}$&$q^{-1/6}$&$-q^{-1/3}$&$-q^{-1/2}$\\
	    \hline
	\end{tabular}
\end{center}
			
	Here \(\mathfrak a\) denotes a cusp of the modular curve associated with \(\Gamma_0(6)\), while \(w_{\mathfrak a}\) denotes the width of that cusp and \(\beta_{\mathfrak a}\) be the positive generator of the
	isotropic lattice after scaling the cusp to infinity.
	
	In the negative-index formula of \cite[Proposition~4.7]{BFtraces},
	an index \(-6t^2\) receives a polar term \(c_{\mathfrak a}(-r)\)
	only when \(r\in(2t/\beta_{\mathfrak a})\Z_{>0}\).
	For a polar exponent \(r=1/w_{\mathfrak a}\), this requires
	\(12t j=1\) for some positive integer \(j\).  Thus \(t\le1/12\).
	In component \(1\), admissible negative indices are congruent to
	\(-1/24\) modulo \(\Z\), so only \(t=1/12\) remains.
	
	At this index the oriented cusp counts are
	\[
	v_{\mathfrak a}(1,-1/24)=
	\begin{cases}1,&\mathfrak a=\infty,\\0,&\text{otherwise},\end{cases}
	\qquad
	v_{\mathfrak a}(-1,-1/24)=
	\begin{cases}1,&\mathfrak a=0,\\0,&\text{otherwise}.\end{cases}
	\]
	The two translation parameters in the exponential factors are zero.
	The negative coefficient is therefore
	$
	-c_\infty(-1)-c_0(-1/6)=-1-1=-2.$
	These are the specialized cusp data in \cite[Section~3]{AA}. This completes the proof of the lemma.
\end{proof}

\subsection{Identification of the lift and proof of the trace formula}
Combining Lemmas~\ref{lem:factor2} and \ref{lem:principal}, we get
\begin{equation}\label{eq:Iliftseries}
	\frac12I_1(\tau,f)^+
	=-q^{-1/24}
	+\sum_{n\ge1}\Tr_{D_n}(f)q^{n-1/24}.
\end{equation}
The harmonic Maass forms \(I_1/2\) and \(\Mcal\) have the same
multiplier and the same principal part.  Their shadows are holomorphic
at the cusp by the theta-lift theorem and
Theorem~\ref{thm:completion}.  Hence, Proposition~\ref{lem:unique} yields
\begin{equation}\label{eq:liftidentity}
	I_1(\tau,f)=2\Mcal(\tau).
\end{equation}
Now by comparing the coefficients on both sides proves \eqref{eq:ftrace}. Next, by subtracting the Bruinier--Ono formula \eqref{eq:partitionCM}
gives \eqref{eq:sptCM}, and hence the identities in
Theorem~\ref{thm:CMmain}.
\subsection{Algebraicity and class polynomials}\label{subsec:algebraicity}
\label{subsec:algebraicity}
We next establish the algebraicity properties of the CM values
appearing in the trace formula. For \(Q=[a,b,c]\in \Qcal_{D_n}^{(1)}\), let
$t=\gcd(a,b,c), Q_0=\frac{Q}{t},
 d=\frac{D_n}{t^2}.
$
Since \(D_n\equiv1\pmod{24}\), we have
\(\gcd(t,6)=1\). In particular, \(t\) is invertible modulo
\(12\). If \(b\equiv1\pmod{12}\), then the middle coefficient
of \(Q_0\) satisfies
$
\frac bt\equiv t^{-1}\pmod{12}.
$
Moreover, \(Q_0\) is primitive and has discriminant \(d\).\\
\noindent
The discriminant \(d\) is the discriminant of the quadratic order
associated with \(Q_0\). Since \(t^2\equiv1\pmod{24}\), it follows
that \(d\equiv1\pmod{24}\). Consequently, the primes \(2\) and
\(3\) split in the underlying imaginary quadratic field. This
splitting condition is compatible with the level-six CM
construction and permits the relevant level structure to be
described in terms of CM points with the prescribed orientation.\\
\noindent
Let \(H_d\) denote the ring class field of the quadratic order of
discriminant \(d\). By the theory of complex multiplication,
together with the level structure and invariance properties of
the modular function \(f\), the value \(f(\alpha_Q)\) belongs to
the appropriate ring class field. The corresponding algebraicity
of \(P(\alpha_Q)\) follows from the CM-value results of
Bruinier--Ono~\cite{BO}. In applying these results, the
modular functions and the CM points must be taken with the same
level structure and orientation.\\
\noindent
The order class group acts on the CM points associated with a
fixed discriminant \(d\), and hence permutes the corresponding
CM values. Since \(d=D_n/t^2\), the order of discriminant \(d\)
has conductor dividing the conductor of the order of
discriminant \(D_n\). The standard inclusions of ring class fields
therefore give
$
H_d\subseteq H_{D_n}.
$
Thus the values arising from each fixed-content family lie in
the ring class field associated with \(D_n\).\\
\noindent
We now briefly discuss integrality. The weakly holomorphic modular
function \(f\) has integral Fourier expansions at all cusps.
The relevant modular-function integrality criterion
\cite[Lemma~4.3]{BO} implies that
$f(\alpha_Q)\in\overline{\mathbb Z},$ the ring of algebraic integers in \(\overline{\mathbb Q}\).
For completeness, the underlying argument is as follows. Form
the modular polynomial
\[
\Psi_f(X,\tau)
:=
\prod_{\gamma\in\Gamma_0(6)\backslash
	\operatorname{SL}_2(\mathbb Z)}
\bigl(X-f|_0\gamma(\tau)\bigr).
\]
Its coefficients are modular functions for
\(\operatorname{SL}_2(\mathbb Z)\), holomorphic on the upper
half-plane. They are therefore polynomials in \(j(\tau)\).
The integral Fourier expansions at all cusps imply that the
coefficients of these polynomials are algebraic integers.
Since \(j(\alpha_Q)\) is an algebraic integer, the value
\(f(\alpha_Q)\) satisfies a monic polynomial over
\(\overline{\mathbb Z}\), and is itself an algebraic integer.

For the nonholomorphic function \(P\), algebraicity alone does
not imply integrality. The stronger denominator theorem of
Larson--Rolen~\cite[Theorem~1.1]{LR} applies to
the relevant weight-zero weak Maass form, provided the required
integrality conditions hold for the Fourier expansions at every
cusp. It gives
$
D_n P(\alpha_Q)\in\overline{\mathbb Z}.
$
Combining this with the integrality of \(f(\alpha_Q)\), we obtain
\begin{equation}\label{eq:gIntegral}
	D_n\bigl(f(\alpha_Q)-P(\alpha_Q)\bigr)
	\in\overline{\mathbb Z}.
\end{equation}
\noindent
Finally, we explain why the values in the trace formula form a
Galois-stable multiset. Complex conjugation reverses the residue \(b\bmod12\); the \(W_6\)
involution reverses it back and fixes both functions.  Together with
CM reciprocity, this proves that the multiset of values in
\eqref{eq:sptCM} is stable under Galois conjugation over \(\Q\).
If
\[
h_n^\ast=\#\bigl(\Gamma\backslash\Qcal_{D_n}^{(1)}\bigr),
\]
then the monic class polynomial
\begin{equation}\label{eq:classpoly}
	\mathscr H_n(X)=
	\prod_{Q\in\Gamma\backslash\Qcal_{D_n}^{(1)}}
	\bigl(X-(f-P)(\alpha_Q)\bigr)
	=
	X^{h_n^\ast}-12\spt(n)X^{h_n^\ast-1}+\cdots
\end{equation}
belongs to \(\Q[X]\).  Scaling each root by \(24n-1\) gives a
monic polynomial in \(\Z[X]\).  This completes the algebraicity
statement of Theorem~\ref{thm:CMmain}.
\begin{example}[Discriminant \(-23\)]\label{ex:minus23}
	For \(n=1\), the three forms may be represented by
	\[
	[6,1,1],\qquad [12,13,4],\qquad [18,25,9].
	\]
	In this case, \cite[Example following Theorem~2]{AA} gives
	\[
	\mathscr H_1(X)=X^3-12X^2-\frac{1008}{23}X-\frac{1728}{23}.
	\]
	The trace of the roots is \(12\), in agreement with \(\spt(1)=1\).
	The corresponding integral polynomial for the scaled roots is
	\[
	\prod_Q\bigl(X-23(f-P)(\alpha_Q)\bigr)
	=X^3-276X^2-23184X-914112.
	\]
\end{example}

\section{Atkin operators and the Fricke transformation}\label{sec:atkin}

Through out this section, let $\ell\geq5$ be prime, and write
\begin{align}\label{esp}
s_\ell:=\frac{\ell^2-1}{24},\qquad
\chi:=\chi_{12}(\ell)=\left(\frac{12}{\ell}\right),\qquad
\epsilon_\ell=\begin{cases}1,&\ell\equiv1\pmod4,\\
	i,&\ell\equiv3\pmod4.\end{cases}
\end{align}
We use the convention $a(x)=0$ whenever $x$ is not a nonnegative integer.
The construction below is Garvan's construction
\cite[Theorem~2.2, Corollary~2.3, and Theorem~4.1]{Garvan}.
We give the matrix and Fourier calculations in order to fix every normalization. Recall that
\begin{equation}\label{fr:completion}
	\mathcal M(\tau)=\alpha(\tau)+\mathcal M^-(\tau),\qquad
	\mathcal M^-(\tau)=\frac3{\sqrt\pi}
	\sum_{m\geq1}\chi_{12}(m)m
	\Gamma\!\left(-\frac12,\frac{\pi m^2y}{6}\right)q^{-m^2/24},
\end{equation}
where $y=\operatorname{Im}\tau$, and that
\begin{align}\label{fr:Mtrans}
	\mathcal M(\gamma\tau)
	=\nu_\eta(\gamma)^{-1}(c\tau+d)^{3/2}\mathcal M(\tau),
	\qquad \gamma=\begin{pmatrix}a&b\\c&d\end{pmatrix}\in
	\operatorname{SL}_2(\mathbb Z).
\end{align}
Fractional powers use the principal branch on the upper half-plane.
In particular,
\begin{equation}\label{fr:Sinversion}
	\mathcal M(-1/\tau)=e^{\pi i/4}\tau^{3/2}\mathcal M(\tau).
\end{equation}
For series with fractional exponents, set
\[
U_\ell^* H(\tau):=\frac1\ell\sum_{j=0}^{\ell-1}
H\!\left(\frac{\tau+24j}{\ell}\right).
\]

\begin{lemma}\label{fr:annihilation}
	The function
	\begin{equation}\label{fr:alpha}
		\alpha_\ell(\tau):=U_\ell^*\mathcal M(\tau)
		-\ell.\chi.\mathcal M(\ell\tau)
		=\sum_{n\geq0}
		\left(a(\ell n-s_\ell)-\ell\chi a(n/\ell)\right)
		q^{n-\ell/24}
	\end{equation}
	is holomorphic on the upper half-plane.
\end{lemma}

\begin{proof} We prove the lemma by separating the holomorphic and
	nonholomorphic parts of \(\mathcal M\). The essential point is that
	the averaging operator \(U_\ell^*\) annihilates all nonholomorphic
	Fourier terms whose indices are not divisible by \(\ell\), and
	the remaining terms cancel against
	\(\ell\chi\mathcal M^-(\ell\tau)\).
	By averaging $q^{-m^2/24}$ over the indicated shifts produces the factor
	$\ell^{-1}\sum_{j=0}^{\ell-1}e^{-2\pi i m^2j/\ell}$.
	This factor is one when $\ell\mid m$ and zero otherwise.
	Put $m=\ell r$ in the surviving terms of
	\eqref{fr:completion}. The identities
	$\chi_{12}(\ell r)=\chi\chi_{12}(r)$ and
	\[
	\Gamma\!\left(-\frac12,
	\frac{\pi(\ell r)^2(y/\ell)}6\right)
	=\Gamma\!\left(-\frac12,\frac{\pi r^2(\ell y)}6\right)
	\]
	then give
	\[
	U_\ell^*\mathcal M^-(\tau)
	=\ell\chi\mathcal M^-(\ell\tau).
	\]
	Consequently only holomorphic terms remain in
	\eqref{fr:alpha}. On its holomorphic part, the averaging condition is
	$24m\equiv1\pmod\ell$. Since
	$24s_\ell=\ell^2-1$, write $m=\ell n-s_\ell$; then
	$(m-1/24)/\ell=n-\ell/24$. The second term in
	\eqref{fr:alpha} contributes exactly when $\ell\mid n$.
	This proves the displayed Fourier expansion.
\end{proof}
\begin{lemma}\label{fr:modularity-interior}
	The function
	\begin{equation}\label{fr:Gdef}
		G_\ell(\tau):=\alpha_\ell(\tau)
		\frac{\eta(\tau)^{2\ell}}{\eta(\ell\tau)}
	\end{equation}
	is holomorphic on the upper half-plane and transforms with weight
	$\ell+1$ and trivial character on $\Gamma_0(\ell)$.
\end{lemma}

\begin{proof}
	Note that the product $\mathcal M\eta$ transforms with weight $2$ on
	$\operatorname{SL}_2(\mathbb Z)$. Using Theorem 1.64 in \cite{CBMS}, we find that the eta quotient
	$\eta(\ell^2\tau)/\eta(\tau)$ transforms with weight zero and
	trivial character on $\Gamma_0(\ell^2)$: the two eta congruence
	conditions are $\ell^2-1\equiv0\pmod{24}$ and
	$1-\ell^2\equiv0\pmod{24}$, and its character is trivial since
	$\ell^2$ is a square. Thus
	\[
	\mathcal F_\ell(\tau)
	:=\mathcal M(\tau)\eta(\ell^2\tau)
	\]
	transforms with weight $2$ on $\Gamma_0(\ell^2)$.
	Moreover,
	\[
	U_\ell^*\mathcal F_\ell
	=\eta(\ell\tau)U_\ell^*\mathcal M,
	\]
	since $\eta(\ell\tau+24\ell j)=\eta(\ell\tau)$.\\
	\noindent
	Next, we show the level-lowering calculation explicitly. Let
	$\gamma=\begin{pmatrix}a&b\\c&d\end{pmatrix}\in\Gamma_0(\ell)$,
	and put $B_j:=\begin{pmatrix}1&24j\\0&\ell\end{pmatrix}$.
	For each $j$ choose $j'$ modulo $\ell$ by
	$24aj'\equiv b+24jd\pmod\ell$. The correspondence $j\mapsto j'$
	is a permutation of the residue classes modulo $\ell$, and
	\[
	B_j\gamma=\gamma_j B_{j'},\qquad
	\gamma_j=
	\begin{pmatrix}
		a+24jc&\dfrac{b+24jd-24j'(a+24jc)}{\ell}\\[4pt]
		\ell c&d-24j'c
	\end{pmatrix}\in\Gamma_0(\ell^2).
	\]
	The factor from the lower row is
	$\ell c B_{j'}\tau+d-24j'c=c\tau+d$.
	Summing the transformations of $\mathcal F_\ell$ therefore shows that
	$\eta(\ell\tau)U_\ell^*\mathcal M$ has weight $2$ on
	$\Gamma_0(\ell)$. The same is true of
	$\mathcal M(\ell\tau)\eta(\ell\tau)$, using
	$\begin{pmatrix}a&\ell b\\c/\ell&d\end{pmatrix}\in
	\operatorname{SL}_2(\mathbb Z)$.
	Their difference gives weight $2$ for
	$\alpha_\ell(\tau)\eta(\ell\tau)$.\\
	\noindent
	Again by using Theorems 1.64 and 1.65 in \cite{CBMS}, we find that
	$
	\frac{\eta(\tau)^{2\ell}}{\eta(\ell\tau)^2}
	$
	has weight $\ell-1$ and trivial character on $\Gamma_0(\ell)$.
	Its character parameter is
	$\ell^{-2}$, a rational square. Multiplying gives the required
	transformation of $G_\ell$. Holomorphy on the upper half-plane follows
	from Lemma~\ref{fr:annihilation} and the nonvanishing of $\eta$ there.
\end{proof}
Next, we define
\begin{equation}\label{fr:betadef}
	\beta_\ell(\tau):=\sum_{n\geq-s_\ell}
	\left\{
	\chi a(n)\left[\left(\frac{1-24n}{\ell}\right)-1\right]
	+\ell a\!\left(\frac{n+s_\ell}{\ell^2}\right)
	\right\}q^{n-1/24},
\end{equation}
and
\begin{equation}\label{fr:Jdef}
	J_\ell(\tau):=\beta_\ell(\tau)
	\frac{\eta(\ell\tau)^{2\ell}}{\eta(\tau)}.
\end{equation}

\begin{proposition}\label{fr:Fricke}
	The following identities and integrality assertions hold:
	\begin{align}
		\alpha_\ell\!\left(-\frac1{\ell\tau}\right)
		&=\ell e^{\pi i/4}\tau^{3/2}\beta_\ell(\tau),
		\label{fr:alphafr}\\
		G_\ell\!\left(-\frac1{\ell\tau}\right)
		&=-(i\ell\tau)^{\ell+1}J_\ell(\tau),
		\label{fr:Gfr}\\
		G_\ell(\tau)&\in\ell\chi+q\mathbb Z[[q]],\qquad
		J_\ell(\tau)=-\ell q^{s_\ell}+O(q^{s_\ell+1})
		\in q^{s_\ell}\mathbb Z[[q]].\label{fr:integrality}
	\end{align}
	In particular, $G_\ell\in M_{\ell+1}(\Gamma_0(\ell))$.
\end{proposition}

\begin{proof}
	Let $S=\begin{pmatrix}0&-1\\1&0\end{pmatrix}$.
	For $1\leq j\leq\ell-1$, choose $j^*$ with
	$576jj^*\equiv-1\pmod\ell$. Then for $B_j=\begin{pmatrix}1&24j\\0&\ell\end{pmatrix}$, we have that
	\begin{equation}\label{fr:matrixS}
		B_jS=C_jB_{j^*},\qquad
		C_j=\begin{pmatrix}
			24j&(-1-576jj^*)/\ell\\
			\ell&-24j^*
		\end{pmatrix}\in\operatorname{SL}_2(\mathbb Z).
	\end{equation}
	For odd $c>0$, the explicit eta multiplier formula is
	\[
	\nu_\eta\!\begin{pmatrix}a&b\\c&d\end{pmatrix}
	=\left(\frac d c\right)
	\exp\!\left(\frac{\pi i}{12}
	\big((a+d)c-bd(c^2-1)-3c\big)\right).
	\]
	For $C_j$, both $a$ and $d$ are divisible by $24$, so it specializes to
	\begin{equation}\label{fr:etamult}
		\nu_\eta(C_j)=
		\left(\frac{-24j^*}{\ell}\right)e^{-\pi i\ell/4}.
	\end{equation}
	Also, $\ell B_{j^*}z-24j^*=z$. Hence
	\eqref{fr:Mtrans} and \eqref{fr:matrixS} give
	\[
	\mathcal M(B_jSz)=z^{3/2}
	\left(\frac{-24j^*}{\ell}\right)e^{\pi i\ell/4}
	\mathcal M(B_{j^*}z).
	\]
	Using this with $z=\ell\tau$ for the nonzero terms in
	$U_\ell^*\mathcal M(-1/(\ell\tau))$. For the term with $j=0$
	and the subtraction term, use \eqref{fr:Sinversion}.
	After factoring out $\ell e^{\pi i/4}\tau^{3/2}$, we obtain that
	\begin{equation}\label{fr:preGauss}
		\alpha_\ell\!\left(-\frac1{\ell\tau}\right)
		=\ell e^{\pi i/4}\tau^{3/2}
		\left\{\ell\mathcal M(\ell^2\tau)
		+\mathcal R_\ell\mathcal M(\tau)-\chi\mathcal M(\tau)\right\},
	\end{equation}
	where
	\[
	\mathcal R_\ell H(\tau):=
	\frac{e^{\pi i(\ell-1)/4}}{\sqrt\ell}
	\sum_{j=1}^{\ell-1}\left(\frac{-24j}{\ell}\right)
	H\!\left(\tau+\frac{24j}{\ell}\right).
	\]
	
	To evaluate this operator, we use the quadratic Gauss sum
	\[
	\sum_{j=1}^{\ell-1}\left(\frac j\ell\right)e^{2\pi i rj/\ell}
	=\epsilon_\ell\sqrt\ell\left(\frac r\ell\right)
	\qquad(r\in\mathbb Z),
	\]
	including $r\equiv0\pmod\ell$, where $\epsilon_\ell$ is defined in \eqref{esp}. The elementary identity
	$e^{\pi i(\ell-1)/4}\epsilon_\ell=
	\left(\frac{2}{\ell}\right),$ yields
	\[
	e^{\pi i(\ell-1)/4}\epsilon_\ell
	\left(\frac{24}{\ell}\right)=\left(\frac{48}{\ell}\right)
	=\left(\frac{12}{\ell}\right)=\chi.
	\]
	Consequently the multiplier on $q^{r/24}$ is
	\begin{equation}\label{fr:Gaussmult}
		\frac{e^{\pi i(\ell-1)/4}}{\sqrt\ell}
		\sum_{j=1}^{\ell-1}\left(\frac{-24j}{\ell}\right)
		e^{2\pi i rj/\ell}
		=\chi\left(\frac{-r}{\ell}\right).
	\end{equation}
	Thus $\mathcal R_\ell$ multiplies the holomorphic coefficient of
	$q^{n-1/24}$ by $\chi((1-24n)/\ell)$.\\
	\noindent 
	Next, we calculate the contributions of nonholomorphic part. In
	\eqref{fr:completion} exponents of $\mathcal M^-(\tau)$ have $r=-m^2$, so
	\eqref{fr:Gaussmult} gives multiplier $\chi$ when $\ell\nmid m$
	and zero when $\ell\mid m$. The former terms cancel against
	$-\chi\mathcal M^-$. Let $\mathcal N_{\mathrm{div}}$ denote the
	subseries of $\mathcal M^-$ restricted to $\ell\mid m$.
	Putting $m=\ell r$ without rescaling $\tau$ gives
	\[
	\mathcal N_{\mathrm{div}}(\tau)
	=\ell\chi\mathcal M^-(\ell^2\tau).
	\]
	Since $\chi^2=1$, the remaining contribution is
	\[
	\ell\mathcal M^-(\ell^2\tau)-\chi
	\mathcal N_{\mathrm{div}}(\tau)=0.
	\]
	The expression in braces in \eqref{fr:preGauss} is therefore
	holomorphic. Its holomorphic Fourier expansion is
	\eqref{fr:betadef}, This proves \eqref{fr:alphafr}.\\
	\noindent
	The eta identities
	\[
	\eta\!\left(-\frac1{\ell\tau}\right)
	=(-i\ell\tau)^{1/2}\eta(\ell\tau),\qquad
	\eta(-1/\tau)=(-i\tau)^{1/2}\eta(\tau)
	\]
	now give \eqref{fr:Gfr}; the scalar simplifies because $\ell$ is odd:
	\[
	\ell e^{\pi i/4}\tau^{3/2}
	\frac{(-i\ell\tau)^\ell}{(-i\tau)^{1/2}}
	=-(i\ell\tau)^{\ell+1}.
	\]
	
	For the expansions, define $\mathscr E(q):=\prod_{m\geq1}(1-q^m)$.
	Then
	\[
	G_\ell(\tau)=
	\frac{\mathscr E(q)^{2\ell}}{\mathscr E(q^\ell)}
	\sum_{n\geq0}\left(a(\ell n-s_\ell)-\ell\chi a(n/\ell)\right)q^n.
	\]
	The coefficient at $n=0$ is $\ell\chi$, since $a(0)=-1$.
	Likewise
	\[
	J_\ell(\tau)=\frac{\mathscr E(q^\ell)^{2\ell}}{\mathscr E(q)}
	\sum_{n\geq-s_\ell}
	\left\{\chi a(n)\left[\left(\frac{1-24n}{\ell}\right)-1\right]
	+\ell a\!\left(\frac{n+s_\ell}{\ell^2}\right)\right\}q^{n+2s_\ell}.
	\]
	Its first coefficient occurs at $n=-s_\ell$ and equals $-\ell$.
	Both product quotients have integral coefficients and constant term one.
	This proves \eqref{fr:integrality}. There are exactly two cusps for
	$\Gamma_0(\ell)$, represented by $\infty$ and $0$. The expansion of
	$G_\ell$ gives holomorphy at $\infty$, and \eqref{fr:Gfr} together with
	the expansion of $J_\ell$ gives holomorphy at $0$. This completes the proof of the lemma.
\end{proof}

\section{The modular trace and proof of the congruences}\label{sec:congruenceproof}

In this section, we prove our main result. For a \(1\)-periodic function with an integral-power Fourier expansion,
we use the ordinary Atkin operator
\[
U_\ell H(\tau)=\frac1\ell\sum_{j=0}^{\ell-1}
H\!\left(\frac{\tau+j}{\ell}\right),\qquad
U_\ell\!\left(\sum_n c(n)q^n\right)=\sum_n c(\ell n)q^n.
\]
The shifted operator \(U_\ell^*\) used in Section~\ref{sec:atkin}
is adapted instead to the fractional exponents of \(\Mcal\). For an integer weight $k$, we recall the slash operator:
\[
(H|_k\gamma)(\tau):=\det(\gamma)^{k/2}(c\tau+d)^{-k}H(\gamma\tau),
\qquad \det(\gamma)>0,
\]
and put $W_\ell=\begin{pmatrix}0&-1\\\ell&0\end{pmatrix}$.
It follows from \eqref{fr:Gfr} that
\begin{equation}\label{fr:slashFricke}
	G_\ell|_{\ell+1}W_\ell
	=(-1)^{(\ell-1)/2}\ell^{(\ell+1)/2}J_\ell.
\end{equation}

\begin{proposition}\label{fr:trace}
	The form
	\begin{equation}\label{fr:Kdef}
		K_\ell:=G_\ell+(-1)^{(\ell-1)/2}\ell U_\ell J_\ell\in M_{\ell+1}(\operatorname{SL}_2(\mathbb Z))\cap\mathbb Z[[q]]
	\end{equation}
	and it satisfies $K_\ell\equiv G_\ell\pmod\ell$, and its constant term
	is $\ell\chi_{12}(\ell)$.
\end{proposition}

\begin{proof}
	We take a coset representatives for $\Gamma_0(\ell)\backslash
	\operatorname{SL}_2(\mathbb Z)$ are
	$I$ and $ST^j$ for $0\leq j<\ell$.
	For $A_j=\begin{pmatrix}1&j\\0&\ell\end{pmatrix}$,
	the identity $W_\ell A_j=\ell ST^j$ and the trivial action of positive
	scalar matrices under the slash operator give
	\[
	\sum_{j=0}^{\ell-1}H|_kST^j
	=\ell^{-k/2}\sum_{j=0}^{\ell-1}
	(H|_kW_\ell)\!\left(\frac{\tau+j}{\ell}\right)
	=\ell^{1-k/2}U_\ell(H|_kW_\ell).
	\]
	Hence, the unnormalized group trace (see, for example,
	Diamond--Shurman~\cite{DiamondShurman2005} for the general theory) is
	\[
	\mathcal T_\ell(H)
	:=\operatorname{Tr}_{\Gamma_0(\ell)}^{\operatorname{SL}_2(\mathbb Z)}H
	=H+\ell^{1-k/2}U_\ell(H|_kW_\ell).
	\]
	It is a holomorphic modular form for the full modular group.
	Apply this to $H=G_\ell$ and $k=\ell+1$, and use
	\eqref{fr:slashFricke}; the result is \eqref{fr:Kdef}.
	In particular, $K_\ell=\mathcal T_\ell(G_\ell)$.
	Integrality and the congruence follow from
	\eqref{fr:integrality}. Since $s_\ell>0$, the series $U_\ell J_\ell$
	has zero constant term, proving the final assertion.
\end{proof}

\begin{proof}[Proof of Theorem~\ref{thm:congruences}]
	The structure theorem for full-group modular forms gives
	$M_*(\operatorname{SL}_2(\mathbb Z))=\mathbb C[E_4,E_6]$
	\cite[Chapter~VII]{Serre}. Hence
	\[
	M_6=\mathbb C E_6,\qquad M_8=\mathbb C E_4^2,
	\qquad M_{14}=\mathbb C E_4^2E_6.
	\]
	These generators have integral coefficients and constant term one.
	The constant term in Proposition~\ref{fr:trace} consequently gives
	\begin{equation}\label{fr:exactK}
		K_5=-5E_6,\qquad K_7=-7E_4^2,\qquad
		K_{13}=13E_4^2E_6.
	\end{equation}
	Thus $G_\ell\equiv0\pmod\ell$ for $\ell=5,7,13$.
	For every index occurring in \eqref{fr:alpha},
	\[
	24(\ell n-s_\ell)-1=\ell(24n-\ell),
	\]
	so the partition-number summand in
	$a(\ell n-s_\ell)=12\operatorname{spt}(\ell n-s_\ell)
	+(24(\ell n-s_\ell)-1)p(\ell n-s_\ell)$ vanishes modulo $\ell$.
	Since $12$ and $\mathscr E(q)^{2\ell}/\mathscr E(q^\ell)$ are units in
	$\mathbb F_\ell[[q]]$, we obtain
	\[
	\operatorname{spt}(\ell n-s_\ell)\equiv0\pmod\ell.
	\]
Here $s_5=1$, $s_7=2$, and $s_{13}=7$; replacing $n$ by $n+1$
	yields the three congruences in Theorem \ref{thm:congruences}.
\end{proof}

\section{Specialization of the CM trace to supersingular points}
\label{sec:geometry}

The algebraic trace in \eqref{eq:ftrace} has a direct interpretation
after reduction of the underlying elliptic curves.  In this section
we make that interpretation precise, including the integrality of the
function being evaluated and the hypothesis on the CM field.

\subsection{The integral level-six function}
Let \(Y_0(6)\) denote the noncuspidal modular curve.  Its geometric
points in characteristic prime to \(6\) classify elliptic curves
equipped with a cyclic subgroup of order \(6\).  A weight-zero modular
function that is regular on this integral model can be evaluated on
such a pair without choosing an invariant differential. Here, we adopt the standard algebraic definition of modular forms as sections of tensor powers of the Hodge bundle, following Katz \cite{Katz}.

\begin{lemma}\label{lem:fregular}
	The function \(f\) of \eqref{eq:f} extends to a regular function on
	\(Y_0(6)_{\Z[1/6]}\).  In particular, it is regular at every
	supersingular level-six point in characteristic \(\ell\ge5\).
\end{lemma}

\begin{proof}
	We begin by working with the normal moduli stack of elliptic curves equipped with a cyclic subgroup of order 6.  Write \(\omega\) for its Hodge bundle.  The eta-product \(H_6\) has trivial character and a
	rational \(q\)-expansion at the rational cusp \(\infty\).
	The \(q\)-expansion principle therefore identifies it with an algebraic
	section of \(\omega^{\otimes4}\) over \(\Q\)
	\cite[Sections~1.5--1.6]{Katz}.  We regard this section as a rational
	section on the integral model.  It satisfies
	\begin{equation}\label{eq:H6discriminants}
		H_6(\tau)^{12}
		=\Dmod(\tau)\Dmod(2\tau)\Dmod(3\tau)\Dmod(6\tau),
		\qquad \text{where}~~\Dmod(\tau)=\eta(\tau)^{24}.
	\end{equation}
	This is an identity of rational sections of
	\(\omega^{\otimes48}\).\\
	\noindent
	For \(d\mid6\), we define \(E_\tau:=\C/(\Z+\tau\Z)\) and the cyclic subgroup
	generated by \(1/d\).  The quotient map, after identifying the target
	with \(E_{d\tau}\), is \(z\mapsto dz\); it pulls back the target
	invariant differential to \(d\,dz\).  Consequently the analytic
	section \(E_k(d\tau)\), expressed in the source differential, is
	\(d^{-k}\) times the pulled-back weight-\(k\) section on the quotient.
	The same statement holds for \(\Dmod(d\tau)\) with \(k=12\).
	All these normalization factors are units over \(\Z[1/6]\).
	The quotient is again a smooth elliptic curve on the noncuspidal
	locus, so each of these discriminant sections is nowhere zero there.\\
	\noindent
	By locally trivializing \(\omega\), and letting \(h\) be the rational function
	representing \(H_6\) in a normal local ring \(R\), equation
	\eqref{eq:H6discriminants} gives \(h^{12}\in R^\times\).
	Both \(h\) and \(h^{-1}\) are therefore integral over \(R\); normality
	implies that \(h\in R^\times\).  Thus \(H_6\) is a unit section
	throughout the noncuspidal locus.\\
	\noindent
	Similarly, each \(E_4(d\tau)\), for \(d\mid6\), extends as a regular
	section of \(\omega^{\otimes4}\), by evaluation of \(E_4\) on the
	quotient elliptic curve.  Since \(24\) is invertible over \(\Z[1/6]\),
	the quotient in \eqref{eq:f} is regular.  Being of weight zero, it
	descends to a regular function on the coarse modular curve. This completes the proof of the lemma.
\end{proof}
\subsection{CM reduction and the multiplicities in the trace}
For fix \(n\ge1\), let \(D:=D_n=1-24n\), and define
$
K:=\Q(\sqrt D),\quad
\mathscr C_D=\Gamma\backslash\Qcal_D^{(1)}.
$
The class of \(Q\) gives a level-six CM point \(x_Q\) whose analytic
coordinate is \(\alpha_Q\).  All these elliptic curves have the same
CM field \(K\), although their endomorphism orders may differ when
imprimitive forms occur.\\
\noindent
Choose a number field over which all the finitely many pairs \(x_Q\)
are defined.  Since CM \(j\)-invariants are algebraic integers, the
curves have potentially good reduction at every finite prime; see
\cite{Silverman}.  After a finite extension, they all have good
reduction at a chosen prime \(\lambda\) above a fixed rational prime
\(\ell\ge5\).  Because \(\ell\nmid6\), their cyclic level-six
subgroups extend and specialize.  Fix an embedding of the residue
field at \(\lambda\) into \(\overline{\F}_\ell\).\\
\noindent
Let \(\Supersing_\ell(6)\) be the finite set of geometric isomorphism classes
of supersingular elliptic curves with a cyclic subgroup of order
\(6\).  If all the \(x_Q\) have supersingular reduction, we define
\begin{equation}\label{eq:actualmultiplicities}
	m_{D,\lambda}(x)
	:=\#\{[Q]\in\mathscr C_D:
	\operatorname{red}_\lambda(x_Q)=x\},
	\qquad x\in\Supersing_\ell(6).
\end{equation}
These are ordinary nonnegative integer multiplicities.  Every class
in the original trace, including every imprimitive class, is counted
once. In characteristic $0$, combining the algebraic trace formula of Bruinier and Ono \cite{BO} for $p(n)$ with the formula of Ahlgren and Andersen \cite{AA} for $\spt(n)$ cancels the non-holomorphic Maass shadow, yielding the exact holomorphic trace identity \eqref{eq:ftrace}:
\[
12\spt(n) + (24n-1)p(n) = \sum_{Q \in \Gamma_0(6)\backslash\mathcal{Q}_D} f(\alpha_Q).
\]
Our next result specializes this identity in positive characteristic. 
\begin{theorem}\label{thm:SSspecialization}
	Let \(n\ge1\), \(D=1-24n\), and \(K=\Q(\sqrt D)\).  Suppose that
	\(\ell\ge5\) is inert or ramified in \(K\).  With the choices above,
	\begin{equation}\label{eq:SSftrace}
		12\spt(n)+(24n-1)p(n)
		\equiv
		\sum_{x\in\Supersing_\ell(6)}m_{D,\lambda}(x)\,\overline f(x)
		\pmod\lambda.
	\end{equation}
	In particular, if \(\ell\mid24n-1\), then
	\begin{equation}\label{eq:SSspttrace}
		12\spt(n)
		\equiv
		\sum_{x\in\Supersing_\ell(6)}m_{D,\lambda}(x)\,\overline f(x)
		\pmod\lambda.
	\end{equation}
\end{theorem}

\begin{proof}
	Deuring's reduction theorem \cite{Deuring} says that a CM elliptic
	curve with good reduction has supersingular reduction precisely when
	the residue characteristic is nonsplit in its CM field.  The criterion
	also applies to nonmaximal CM orders.  To reduce that assertion to the
	maximal-order case, take an isogenous curve with multiplication by the
	maximal order of \(K\).  After a finite extension the isogeny and both
	good reductions are defined, and the isogeny specializes.
	Supersingularity is invariant under isogeny, so the criterion depends
	on \(K\), rather than on the conductor of the chosen CM order. Thus every \(x_Q\) reduces to a point of \(\Supersing_\ell(6)\).
	By Lemma~\ref{lem:fregular}, we have that \(f(\alpha_Q)\) is
	\(\lambda\)-integral and that
	\[
	f(\alpha_Q)\bmod\lambda
	=\overline f\bigl(\operatorname{red}_\lambda(x_Q)\bigr).
	\]
	Reducing \eqref{eq:ftrace} and grouping the terms according to their
	reductions proves \eqref{eq:SSftrace}.  If \(\ell\mid24n-1\), the
	partition term on the left vanishes modulo \(\lambda\), giving
	\eqref{eq:SSspttrace}. This completes the proof of the theorem.
\end{proof}

\begin{remark}
	Let \(B_\ell(6)\) be the free abelian group on the set
	\(\Supersing_\ell(6)\).  The right side of \eqref{eq:SSftrace} is the
	evaluation pairing between
	\[
	\sum_xm_{D,\lambda}(x)[x]\in B_\ell(6)
	\quad\text{and}\quad
	\bigl(\overline f(x)\bigr)_{x\in\Supersing_\ell(6)}.
	\]
	This description uses the actual reduction map and requires no
	identification of the multiplicities with optimal-embedding counts.
	An explicit formula in terms of quaternionic orders would require
	additional local analysis, especially at primes dividing the
	conductors of the CM orders.  No such formula is used here.
	Theorem~\ref{thm:SSspecialization} is a specialization of the exact
	trace identity; the vanishing of the modular traces in
	\eqref{eq:introK} proves the full arithmetic progressions in
	Theorem~\ref{thm:congruences}.
\end{remark}

\subsection{Why divisibility of the discriminant is insufficient}
Write \(D=f_D^2d_K\), where \(d_K\) is the fundamental discriminant
of \(K\).  A prime dividing \(D\) may divide the conductor \(f_D\)
while splitting in \(K\).  Consequently the condition
\(\ell\mid24n-1\) does not imply the nonsplitting hypothesis of
Theorem~\ref{thm:SSspecialization}.

\begin{example}\label{ex:ordinaryconductor}
	For \(n=74\),
	\[
	D_{74}=-1775=-71\cdot5^2,
	\qquad
	\leg{-71}{5}=\leg{4}{5}=1.
	\]
	Thus \(5\) splits in \(\Q(\sqrt{-71})\), and the CM curves in this
	trace have ordinary reduction after achieving good reduction.
	Nevertheless \(74=5\cdot14+4\), so
	Theorem~\ref{thm:congruences} gives \(\spt(74)\equiv0\pmod5\).
	This example shows why a supersingular argument with only the
	hypothesis \(5\mid24n-1\) would leave part of the progression
	untreated.  The modular trace proof includes these conductor cases.
\end{example}

\subsection{A direct geometric consequence in characteristic five}
The Hasse invariant of elliptic curves in characteristic \(\ell\) is a
modular form of weight \(\ell-1\), with \(q\)-expansion \(1\), whose
zeros are precisely the supersingular points.  We use this standard
description from \cite{Katz,Silverman}.

\begin{proposition}\label{prop:fmod5}
	The reduction of \(f\) vanishes at every point of \(\Supersing_5(6)\).
	Consequently, if \(n\ge1\), \(5\mid24n-1\), and \(5\) is nonsplit
	in \(\Q(\sqrt{1-24n})\), then \(\spt(n)\equiv0\pmod5\).
\end{proposition}

\begin{proof}
	The reduction of \(E_4\) modulo \(5\) is the Hasse invariant: it has
	weight \(4\) and its \(q\)-expansion is \(1\) modulo \(5\), so the
	claim follows from the \(q\)-expansion principle.  If \((E,C)\) is
	supersingular, then all its quotients by subgroups of orders
	\(d\mid6\) are supersingular.  Therefore every section
	\(E_4(d\tau)\) in the numerator of \(f\) vanishes at \((E,C)\).
	The denominator is a unit by Lemma~\ref{lem:fregular}, proving
	\(\overline f(E,C)=0\).  By applying \eqref{eq:SSspttrace} and using the
	invertibility of \(12\) modulo \(5\), we obtain the final assertion.
\end{proof}

The nonsplitting assumption in Proposition~\ref{prop:fmod5} remains
essential for this particular proof.  It is absent from the full
congruence \eqref{eq:spt5}, already established by modular forms.

\section{Comparison with the partition setting}
\label{sec:comparison}

The two main formulas have the same discriminant \(D_n=1-24n\) and
the same level-six CM points.  Their precise relationship is
\begin{equation}\label{eq:comparisontraces}
	\Tr_{D_n}(P)=(24n-1)p(n),\qquad
	\Tr_{D_n}(f-P)=12\spt(n).
\end{equation}
The common CM geometry is useful, but certain analytic and
characteristic-dependent features must be distinguished.

\subsection{The Laplacian of the smallest-parts trace function}
With the convention
\[
\Delta_0=-y^2\left(\frac{\partial^2}{\partial x^2}
+\frac{\partial^2}{\partial y^2}\right),
\]
the Bruinier--Ono function \cite{BO} satisfies \(\Delta_0P=-2P\).  Because \(f\) is holomorphic on \(\HH\),
\(\Delta_0f=0\).  Thus the function \(g=f-P\) satisfies
\begin{equation}\label{eq:glaplacian}
	\Delta_0g=2P,
	\qquad
	\Delta_0(\Delta_0+2)g=0.
\end{equation}
In particular, \(g\) is not an eigenfunction with eigenvalue \(-2\):
that equation together with \eqref{eq:glaplacian} would imply
\(f=0\), contrary to its principal part.  Its modular invariance and
algebraic CM values suffice for \eqref{eq:sptCM}, while any argument
depending specifically on the eigenvalue of \(P\) requires a separate
verification for \(g\).

\subsection{Supersingular elliptic curves in characteristic thirteen}
The distinctive role of \(13\) in the smallest-parts congruences
cannot be inferred by replacing a prime in an argument based solely
on the special invariants \(j=0\) and \(j=1728\).  The relevant
calculation is elementary.

\begin{proposition}\label{prop:SS13}
	The unique supersingular \(j\)-invariant in characteristic \(13\) is
	\(5\).  In particular, it is different from both \(0\) and
	\(1728\bmod13\).
\end{proposition}

\begin{proof}
	Over \(\overline{\F}_{13}\), write an elliptic curve in the form
	\[
	E:y^2=x^3+Ax+B,
	\qquad 4A^3+27B^2\ne0.
	\]
	The Hasse-invariant criterion \cite[Chapter~V]{Silverman} says that
	\(E\) is supersingular if and only if the coefficient of \(x^{12}\)
	in \((x^3+Ax+B)^6\) vanishes.  Only the choices
	\((x^3)^3(Ax)^3\) and \((x^3)^4B^2\) contribute.  Hence this
	coefficient is
	\[
	\frac{6!}{3!3!}A^3+\frac{6!}{4!2!}B^2
	=20A^3+15B^2
	=7A^3+2B^2
	\quad\text{in }\overline{\F}_{13}.
	\]
	Its vanishing is equivalent to \(B^2=3A^3\).  In that case
	nonsingularity gives \(A\ne0\), and
	\[
	j(E)=1728\frac{4A^3}{4A^3+27B^2}
	=12\frac4{4+3}=5
	\quad\text{in }\F_{13}.
	\]
	Conversely, if \(j(E)=5\), then
	\(5(4A^3+B^2)=48A^3\) in characteristic \(13\), whence
	\(B^2=3A^3\).  The Hasse invariant therefore vanishes.  Such curves
	exist, for instance by choosing \(A=1\) and \(B^2=3\) over the
	algebraic closure; the discriminant is nonzero.  Finally,
	\(1728\equiv12\pmod{13}\), so \(5\) is neither exceptional
	invariant.
\end{proof}
\section{Concluding remark}
\label{app:content}

The lack of a primitivity restriction in \eqref{eq:Qset} is essential.
We make its effect on the trace explicit.  For a negative discriminant
\(d\equiv1\pmod{24}\) and a residue \(r\in(\Z/12\Z)^\times\), let
\[
\Qcal_d^{(r),\mathrm{prim}}
:=
\{[a,b,c]:a>0,\ b^2-4ac=d,\ 6\mid a,\ b\equiv r\pmod{12},
\ \gcd(a,b,c)=1\}.
\]
Define
\[
\Tr_d^{(r),\mathrm{prim}}(u)
:=
\sum_{Q\in\Gamma\backslash\Qcal_d^{(r),\mathrm{prim}}}
u(\alpha_Q),\qquad
\Tr_d^{\mathrm{prim}}(u)=\Tr_d^{(1),\mathrm{prim}}(u).
\]
Since $(\Gamma_0(6)/{\pm I})$ is torsion-free, no stabilizer weights are required for the functions considered here.

\begin{proposition}\label{prop:content}
	Let \(D<0\) satisfy \(D\equiv1\pmod{24}\), and let \(u\) be any one of
	\(f\), \(P\), or \(f-P\).  Then
	\begin{equation}\label{eq:contenttrace}
		\Tr_D(u)=
		\sum_{\substack{t\ge1\\t^2\mid D}}
		\ch(t)\Tr_{D/t^2}^{\mathrm{prim}}(u).
	\end{equation}
\end{proposition}
\begin{proof}
	The content \(t=\gcd(a,b,c)\) is invariant under integral changes of
	variables of determinant one.  Since \(D\equiv1\pmod{24}\), every
	\(t\) with \(t^2\mid D\) is prime to \(6\).  In particular,
	\(t^2\equiv1\pmod{24}\), and \(d=D/t^2\equiv1\pmod{24}\).
	
	Dividing a form of content \(t\) by \(t\) induces a bijection of
	\(\Gamma\)-orbits with
	\(\Qcal_d^{(t^{-1}),\mathrm{prim}}\).  Its CM point does not change:
	\[
	\frac{-b+\sqrt D}{2a}
	=\frac{-b/t+\sqrt{D/t^2}}{2a/t}.
	\]
	Because \(t^{-1}\equiv t\pmod{12}\), the required residue is \(t\).
	The Atkin--Lehner bijections and their eigenvalues on \(u\) are
	\[
	\begin{array}{c|cccc}
		r&1&5&7&11\\ \hline
		\text{operator from residue }1&I&W_3&W_2&W_6\\
		\text{eigenvalue on }u&1&-1&-1&1
	\end{array}
	\]
	by \eqref{eq:ALf} and the corresponding identities for \(P\).
	These operators preserve primitivity.  They give
	\[
	\Tr_d^{(r),\mathrm{prim}}(u)
	=\ch(r)\Tr_d^{\mathrm{prim}}(u).
	\]
	Summing over all possible contents proves \eqref{eq:contenttrace}.
\end{proof}

\begin{example}\label{ex:content575}
	Take \(n=24\), so
	$
	D_n=-575=25(-23).
	$
	For example, \([30,25,10]\) belongs to \(\Qcal_{-575}^{(1)}\) and has
	content \(5\).  The only possible contents are \(1\) and \(5\), and
	\(\ch(5)=-1\).  Since \(-23\) has no imprimitive forms,
	\eqref{eq:contenttrace} and Example~\ref{ex:minus23} imply
	\[
	\Tr_{-575}(f-P)
	=\Tr_{-575}^{\mathrm{prim}}(f-P)-12.
	\]
	Thus replacing the trace in Theorem~\ref{thm:CMmain} by a
	primitive-only trace would give \(\spt(24)+1\).
	Likewise,
	\[
	\Tr_{-575}(P)=\Tr_{-575}^{\mathrm{prim}}(P)-23,
	\]
	so dividing the primitive-only trace of \(P\) by \(575\) would give
	\(p(24)+1/25\).  Each fixed-content packet has its own order class-group
	description.  Their union is the trace required by the coefficient
	formula.
\end{example}
\noindent The content decomposition above motivates the following problem of understanding primitive CM traces and their arithmetic properties.
\begin{problem}
	\label{prob:primitive-traces}
	Find explicit formulas for \(\Tr_D^{\mathrm{prim}}(u)\), where
	\(u\in\{f,P,f-P\}\), in terms of class numbers, modular-form
	coefficients, or special values of \(L\)-functions. In particular,
	determine whether the content decomposition
	\[
	\Tr_D(u)=
	\sum_{t^2\mid D}\ch(t)\Tr_{D/t^2}^{\mathrm{prim}}(u)
	\]
	yields new congruences for primitive traces or a uniform
	class-group interpretation of the fixed-content contributions.
\end{problem}

\end{document}